\documentclass[11pt]{article}
\usepackage{amsmath}
\usepackage{amscd}
\usepackage{amsfonts}
\usepackage{amsmath,amsthm,hyperref}
\usepackage[all]{xy}

\newtheorem{theorem}{Theorem}[section]

\newtheorem{definition}[theorem]{Definition}

\theoremstyle{remark}
\newtheorem{remark}[theorem]{Remark}

\begin{document}
\title{\bf{Blaschke's problem for timelike surfaces in pseudo-Riemannian space forms}}
\author{{Peng Wang }
\footnote{Department of Mathematics, Tongji University, 200092,
Shanghai, People's Republic of China.  $~E-mail\ address:
netwangpeng@tongji.edu.cn$ }\footnote{  Supported by Program for Young Excellent Talents in Tongji
University and the Tianyuan Foundation of China (Grant No.10926112).}}
\date{}
\maketitle
\begin{center}
{\bf Abstract}
\end{center}

\noindent{\ \ \small\hspace{2.8mm} We show that isothermic surfaces
and S-Willmore surfaces are also the solutions to the corresponding
Blaschke's problem for both spacelike and timelike surfaces in
pseudo-Riemannian space forms. For timelike surfaces both Willmore
and isothermic, we obtain a description by minimal surfaces similar
to the classical results of Thomsen. 

\section{Introduction}

 Isothermic surfaces and Willmore surfaces are important objects in
 conformal differential geometry. They are both surface classes
invariant under conformal transforms. Although seemed so distinct to
each other, they may be introduced as the only non-trivial solutions
to a problem in the category of conformal differential
geometry\cite{Blaschke}, i.e. the Blaschke's problem:\\

 \textbf{Blaschke's Problem:}
\begin{itshape} Let $S$ be a sphere congruence with
two envelops $$f,f':\ M^2\rightarrow S^3,$$ such that these envelops
induce the same conformal structure. Characterize such sphere
congruences and envelop surfaces.\end{itshape}\\

Blaschke posed this question and solved it in \cite{Blaschke},
namely\\

\textbf{\bf Theorem}
 \begin{itshape} The non-trivial solution to the Blaschke's problem
is either a pair of isothermic surfaces forming Darboux transform to
each other, or a pair of dual Willmore surfaces with their common
mean curvature spheres. (Here non-trivial means the two envelops are
not congruent up to M\"{o}bius transforms.)
\end{itshape}\\

In \cite{Ma2005}, Ma considered the arbitrary co-dimensional case
and proved that the generalized Darboux pair of isothermic surfaces
as well as S-Willmore surfaces in $S^n$ are the full non-trivial
solution to the Blaschke's problem. There have been some other kinds
of generalization as to Blaschke's problem after Ma's work, see
\cite{Dajczer-Tojeiro2007}.

Besides these, isothermic surfaces and Willmore surfaces have close
relations with integrable system. Isothermic surfaces allow many
kinds of transform, which indicate that there is a structure of
integrable system underlying the theory about isothermic surfaces,
see \cite{Bur},\cite{Ma2005}. This fact was just revealed in the
past 20 years. The transforms of Willmore surfaces are more
complicated, but they also can be interpreted by use of integrable
system. Recently, there have been several kinds of research
concerning the surface theory in pseudo-Riemannian space forms from
different viewpoints, for example, see \cite{Deng-Wang},
\cite{Fu-In2003}, \cite{ZCC2004}, \cite{ZCC20042}. For such surfaces
therein, there is also a parallel theory of conformal geometry. When
dealing with the conformal geometry of such surfaces, it is natural
to consider the Blaschke's problem, which is just the main content
of this paper. We obtain the similar results for both spacelike and
timelike surfaces in pseudo-Riemannian space forms as \cite{Ma2005}.
Here $spacelike$ means the induced metric on the surface is
Riemannian and $timelike$ means the induced metric on the surface is
Lorentzian.  For the spacelike case, the theorems and proofs are the
same as Ma's results. So we omit it and just focus on the timelike
case.

As we know that 2-spheres are basic invariants in the conformal
geometry, timelike 2-spheres are the corresponding notions in the
pseudo-Riemannian case. Such 2-spheres comes from the notion of
curvature sphere, which are the 2-spheres tangent to the surface and
have the same mean curvature as the surface in the corresponding
point. For our case, a timelike 2-sphere in $R^3_1$ is just a
hyperboloid of one sheet. Just as the spacelike 2-spheres can be
identified with a 4-dim Lorentzian subspaces in some
$\mathbb{R}^{n+2}_{1}$, timelike 2-sphere can be identified with a
4-dim (2,2)-type subspaces of $\mathbb{R}^{n+2}_{r+1}$. For the
equivalence of them, we refer to \cite{Her1}. So Blaschke's problem
for
timelike surfaces can be stated as below:\\

\textbf{\bf Blaschke's Problem for timelike surfaces:}
 \begin{itshape} Let
$S$ be a timelike 2-sphere congruence with two timelike envelops
$f,f': M^2\rightarrow Q^n_r$, such that these timelike envelops
induce the same conformal structure. Characterize such sphere
congruences and envelop surfaces.
\end{itshape}\\

Here a timelike 2-sphere congruence in $Q^n_r$ is just a map into
the Grassmannian manifold $G_{2,2}(\mathbb{R}^{n+2}_{r+1}):=
\{\text{4-dim (2,2)-type subspaces of}~\mathbb{R}^{n+2}_{r+1}\}$.
For the notion of $Q^n_r$, see Section~2. Then our main result is
the following
theorem:\\

\textbf{Theorem A}
\begin{itshape}
The non-trivial solution to the Blaschke's problem of timelike
surfaces is either a pair of timelike isothermic surfaces forming
Darboux transform to each other, or a pair of dual timelike
S-Willmore surfaces with their common mean curvature spheres.
\end{itshape}\\

We also give a characteristic of timelike isothermic Willmore
surfaces in $Q^3_1$ as follows, which is similar to the classical
results of Thomsen \cite{Thomsen}:\\

\textbf{Theorem B}
\begin{itshape}
Any timelike isothermic Willmore surface in $Q^3_1$ is conformally
equivalent to a timelike minimal surface in some 3-dimensional
Lorentzian space form $R^3_1$, $S^3_1$ or $H^3_1$.
\end{itshape}\\

This paper is organized as follows. In Section~2, the
pseudo-Riemannian conformal space $Q^{n}_{r}$, the general theory
about timelike surfaces in $Q^{n}_{r}$ and the characterization of
timelike Willmore surfaces and timelike isothermic surfaces are
introduced. Then we prove Theorem A, B in Section~3 and Section~4
respectively.

\begin{remark}We note that
there is a parallel and independent paper \cite{Magid} on this topic, which also cites our paper. \cite{Magid} is to deal with timelike surfaces in 3-dimensional Lorenzian space forms. Here our goal is to discuss the interest relationship of isothermic surfaces and S-Willmore surfaces for timelike surfaces in n-dimensional pseudo-Riemannian space forms, following the treatment for spacelike surfaces in \cite{BPP}, \cite{Ma2005}.  \end{remark}

\section{Timelike surfaces in pseudo-Riemannian space forms}

\subsection{conformal geometry of $Q^{n}_{r}$}

Let $\mathbb{R}^m_s$ be the space $\mathbb{R}^m$ equipped with the
quadric form
\[
\langle x,x\rangle=\sum^{m-s}_{1}x^2_i-\sum^m_{m-s+1}x^{2}_i.
\]
We denote by $C^{m-1}_s$ the light cone of $\mathbb{R}^m_s$. The
quadric
\[
Q^n_r=\{\ [x]\in\mathbb{R}P^{n+1}\ |\ x\in C^{n+1}_r\setminus \{0\}
\}
\]
is exactly the projectived light cone. The standard projection
$\pi:C^{n}_{r}\setminus\{0\}\rightarrow Q^n_r$ is a fiber bundle
with fiber $\mathbb{R}\setminus\{0\}$. It is easy to see that
$Q^n_r$ is equipped with a $(n-r,r)-$type pseudo-Riemannian metric
induced from projection $S^{n-r}\times S^r\rightarrow Q^n_r$. Here
\begin{equation}S^{n-r}\times
S^r=\{x\in\mathbb{R}^{n+2}_{r+1}\ |\
\sum^{n-r+1}_{i=1}x^{2}_{i}=\sum^{n+2}_{i=n-r+2}x^2_{i}=1\}\subset
C^{n+1}_r\setminus\{0\}
\end{equation}
endowed with a $(n-r,r)-$type pseudo-Riemannian metric
$g(S^{n-r})\oplus (-g(S^r))$, where $g(S^{n-r})$ and $g(S^r)$ are
standard metrics on $S^{n-r}$ and $S^r$. So there is a conformal
structure of $(n-r,r)-$type pseudo-Riemannian metric $[h]$ on
$Q^n_r$. The conformal group of $(Q^n_r,[h])$ is exactly the
orthogonal group $O(n-r+1,r+1)/\{\pm1\}$, which keeps the inner
product of $\mathbb{R}^{n+2}_{r+1}$ invariant and acts on $Q^n_r$ by
\begin{equation}
T([x])=[xT],\ T\in O(n-r+1,r+1).
\end{equation}

For the three n-dimensional $(n-r,r)-$type pseudo-Riemannian space
forms with constant sectional curvature $c=0,+1,-1$, they are
defined by
\begin{equation*}
R^n_r,\ \  S^n_r :=\{x\in\mathbb{R}^{n+1}_{r}\ |\ \langle
x,x\rangle=1 \},\ \  H^n_r :=\{x\in\mathbb{R}^{n+1}_{r+1}\ |\
\langle x,x\rangle=-1 \}.
\end{equation*}
Each of them could be embedded as a proper subset of $Q^n_r$:
\begin{equation}\label{varphi}
\begin{array}{llll}
\varphi_{0}:R^n_r\rightarrow Q^n_r, ~&
\varphi_{0}(x)=[(\frac{-1+\langle x,x\rangle}{2},
x,\frac{1+\langle x,x\rangle}{2})]; \\[1mm]
\varphi_{+}:S^n_r\rightarrow Q^n_r,
~& \varphi_{+}(x)=[(x,1)]; \\[1mm]
\varphi_{-}:H^n_r\rightarrow Q^n_r, ~& \varphi_{-}(x)=[(1,x)].
\end{array}
\end{equation}
These maps are conformal embeddings. Thus $Q^n_r$ is the proper
space to study the conformal geometry of these $(n-r,r)-$type
pseudo-Riemannian space forms.

\subsection{Basic equations of Timelike Surfaces in $Q^n_r$}

Let $y:M\rightarrow Q^{n}_{r}$ be a timelike surface. For any open
subset $U\subset M$, we call $Y:U\rightarrow C^{n+1}_{r}$ a local
lift of $y$ if $y=\pi\circ Y$. Two different local lifts differ by a
scaling, so the metric induced from them are conformal to each
other. Choose asymptotic coordinates $(u,v)$ on $U$, such that for
some lift $Y$
\begin{equation}\langle Y_{u},Y_{u}\rangle=\langle
Y_{v},Y_{v}\rangle=0.\end{equation} Such property obversely holds
for any lift, showing that asymptotic coordinates are conformal
invariant.

For such a surface there is a decomposition $M\times
\mathbb{R}^{6}_{2}=V\oplus V^{\perp}$, where\begin{equation}
V:=Span\{Y,Y_u,Y_v,Y_{uv}\}
\end{equation}
is a 4-dimensional (2,2)-type subbundle independent of choice of $y$
and $(u,v)$. $V^{\perp}$ has a $(n-r-1,r-1)-$type metric. which
might be identified with the normal bundle of $y$.

For a local asymptotic coordinate $(u,v)$, there is a local lift $Y$
such that $\langle Y_{u},Y_{v}\rangle=\pm\frac{1}{2}$. We can adjust
$v$ so that $$\langle Y_{u},Y_{v}\rangle=\frac{1}{2}.$$ We call $Y$
a canonical lift with respect to $(u,v)$. So there is a unique $N\in
\Gamma(V)$ satisfying
\begin{equation}\label{eq-N-time}
\langle N,Y_{u}\rangle=\langle N,Y_{v}\rangle=\langle
N,N\rangle=0,\langle N,Y\rangle=-1.
\end{equation}

Given frames as above, we note that $Y_{uu}$ and $Y_{vv}$ are
orthogonal to $Y$, $Y_{u}$ and $Y_{v}$. So there must be two
functions $s_1,\ s_2$ and two section $\kappa_1, \kappa_2\in
\Gamma(V^{\perp})$ such that
\begin{equation}\label{eq-Hopf-2}
\left\{\begin {array}{ll}
Y_{uu}=-\frac{s_1}{2}Y+\kappa_{1},\\
Y_{vv}=-\frac{s_2}{2}Y+\kappa_{2}.
\end {array}\right.\end{equation}

These define four basic invariants $\kappa_1$, $\kappa_2$ and $s_1$,
and $s_2$ dependent on $(u,v)$. Similar to the case in M\"obius
geometry, $\kappa_i,$ and $s_i$ are called the \emph{conformal Hopf
differential} and the \emph{Schwarzian derivative} of $y$,
respectively ( compare \cite{BPP},\cite{Ma}).

Let $\psi\in\Gamma(V^{\perp})$ denote a section of the normal
bundle, and $D$ the normal connection, we can derive the structure
equations as below:
\begin{equation}\label{eq-moving2}
\left\{\begin {array}{lllll}
Y_{uu}=-\frac{s_1}{2}Y+\kappa_{1},\\
Y_{vv}=-\frac{s_2}{2}Y+\kappa_{2},\\
Y_{uv}=-\langle \kappa_1,\kappa_2\rangle Y+\frac{1}{2}N,\\
N_{u}=-2\langle \kappa_1,\kappa_2\rangle  Y_{u}-s_{1}Y_{v}+2D_{v}\kappa_1,\\
N_{v}=-s_2Y_{u}-2\langle \kappa_1,\kappa_2\rangle Y_{v}+2D_{u}\kappa_2,\\
\psi_{u}=D_{u}\psi+2\langle\psi,D_{v}\kappa_1\rangle Y-2\langle
\psi,\kappa_1\rangle  Y_{v},\\
\psi_{v}=D_{v}\psi+2\langle\psi,D_{u}\kappa_2\rangle Y-2\langle
\psi,\kappa_2\rangle  Y_{u}.
\end {array}\right.
\end{equation}
The conformal Gauss equations, Codazzi equations, and Ricci
equations as integrable conditions are:
\begin{equation}\label{eq-G2}
\left\{\begin {array}{ll}
\frac{1}{2}s_{1v}=3\langle\kappa_1,D_{u}\kappa_2\rangle+\langle D_{u}\kappa_1,\kappa_2\rangle,\\
\frac{1}{2}s_{2u}=\langle\kappa_1,D_{v}\kappa_2\rangle+3\langle
D_{v}\kappa_1,\kappa_2\rangle;
\end {array}\right.
\end{equation}
\begin{equation}\label{eq-C2}
D_{v}D_{v}\kappa_1+\frac{s_2}{2}\kappa_1=D_{u}D_{u}\kappa_2+\frac{s_1}{2}\kappa_2;
\end{equation}\begin{equation}\label{eq-R2}
R^{D}_{uv}:=D_{u}D_{v}\psi-D_{v}D_{u}\psi =2\langle
\psi,\kappa_1\rangle\kappa_2-2\langle \psi,\kappa_2\rangle\kappa_1.
\end{equation}

\subsection{Timelike Willmore surfaces}

\begin{definition}
Let $y:M\rightarrow Q^{n}_{r}$ be an immersed timelike surface.
\emph{The Willmore functional} of $y$ is defined as:
\[
W(y):=2\int_M\langle\kappa_1,\kappa_2\rangle dudv.
\]
 $y$ is called a
\emph{Willmore surface}, if it is a critical surface of the Willmore
functional with respect to any timelike variation of the map
$y:M\rightarrow Q^{n}_{r}$.
\end{definition}

It is direct to check that $W(y)$ is well-defined. Timelike Willmore
surfaces can be characterized as follows, which is similar to the
results of spacelike case
\cite{AP1996,Deng-Wang,Ma-W-time-Willmore,Wang1998}.

\begin{theorem}\label{thm-willmore-2}
For a timelike surface $y:M^{2}\rightarrow Q^{n}_{r}$, the following
three conditions are equivalent:

(i) $y$ is a timelike Willmore surface.

(ii) The conformal Gauss map $$G:M\rightarrow
G_{2,2}(\mathbb{R}^{n+2}_{r+1}), \ G(p):= V_{p}, \ \forall p\in M$$
of $y$ is harmonic.

(iii) The two Hopf differential $\kappa_1,\kappa_2$ satisfy the
following Willmore equation:
\begin{equation}\label{eq-willmore-time}
D_{v}D_{v}\kappa_1+\frac{s_2}{2}\kappa_1=D_{u}D_{u}\kappa_2+\frac{s_1}{2}\kappa_2=0.
\end{equation}
\end{theorem}

For the proof, we prefer to \cite{Deng-Wang,Ma-W-time-Willmore,
Wang1998}. We also note that the calculation of Euler-Lagrange
equations of Willmore functional by Wang in \cite{Wang1998} is valid
for timelike submanifolds in Lorentzian space forms, and then leads
to Theorem 2.2.

Now we define timelike S-Willmore surfaces as:
\begin{definition}
A timelike Willmore surface $y:M\rightarrow Q^{n}_{r}$ is called an
\emph{S-Willmore surface} if it satisfies $D_{v}\kappa_1\,\|\,
\kappa_1,$ $D_{u}\kappa_2\,\|\, \kappa_2$, i.e., if there exist two
functions $\mu_{1},\mu_{2}$ such that
\begin{equation}D_{v}\kappa_1+\mu_1\kappa_1=D_{u}\kappa_2+\mu_{2}\kappa_2=0.\end{equation}
\end{definition}

\begin{remark}
The notion of S-Willmore comes from Ejiri's work \cite{Ejiri1988}.
He considered the duality theorem for Willmore surfaces in $S^n$,
and found that only some special kinds of Willmore surfaces in
higher dimensional space forms allow dual surfaces. Such surfaces
were named S-Willmore surfaces. Geometrically speaking, the mean
curvature sphere of a S-Willmore surface has two envelop surfaces
which are conformal to each other with opposite orientation. One is
itself and the other is just its dual surface. For the detail, we
refer to \cite{Ejiri1988} and \cite{Her1}.
\end{remark}

\subsection{Timelike isothermic surfaces}

\begin{definition}
Let $y:M\rightarrow Q^{n}_{r}$ be a conformal timelike surface
without umbilic points. It is called $(\pm)-$isothermic if around
each point of $M$ there exists an asymptotic coordinate $(u,v)$ and
canonical lift $Y$ such that the Hopf differentials
$\kappa_1=\pm\kappa_2$. Such a coordinate $(u,v)$ is called an
adapted coordinate.
\end{definition}
$\kappa_1=\pm\kappa_2$ together with the conformal Ricci equations
in \eqref{eq-R2} shows that the normal bundle of $y$ is flat. This
is an important property of isothermic surfaces, which guarantees
that all shape operators commute and the curvature lines could still
be defined. Setting $u=s+t,v=s-t$, the two fundamental forms of an
isothermic surface, with respect to some parallel normal frame
$\{e_{\alpha}\}$, are of the form
\begin{equation}\label{iso-time+}
I=e^{2\rho}(ds^2-dt^2),\ II=\sum_{\alpha}(b_{\alpha1}ds^2-
b_{\alpha2}dt^2)e_{\alpha},\end{equation} if $y$ is $(+)-$isothermic
and
\begin{equation}\label{iso-time-}I=e^{2\rho}(ds^2-dt^2),\ II=\sum_{\alpha}(b_{\alpha1}(ds^2-dt^2)
-b_{\alpha2}dsdt)e_{\alpha}
\end{equation}
if $y$ is $(-)-$isothermic. Note that $(\pm)-$isothermic surfaces
are called real and complex isothermic surface separately in
\cite{Dussan-M2005}. And our notions here follow \cite{Fu-In2003}.

\section{Proof of Theorem A}

Denote $y$ and $\hat{y}$ the pair of surfaces in the Blaschke's
problem. Let $(u,v)$ be an asymptotic coordinate of y and $Y$ the
relevant canonical lift. Choose a lift $\hat{Y}$ of $\hat{y}$ such
that $\langle Y,\hat{Y}\rangle=-1$. Then the sphere congruence
tangent to $Y$ and passing $\hat{Y}$ is
$$Span\{Y,Y_{u},Y_{v},\hat{Y}\}.$$ By the conditions of
Theorem A, we know that $(u,v)$ is also asymptotic coordinate of
$\hat{Y}$, and
\begin{equation}\label{eq-blas}
Span\{Y,Y_{u},Y_{v},\hat{Y}\}=Span\{\hat{Y},\hat{Y}_{u},\hat{Y}_{v},Y\}.\end{equation}
Assume that
\begin{equation}\hat{Y}=N+2aY_u+2bY_v+(2ab+\frac{1}{2}\langle
\xi,\xi\rangle)Y+\xi,
\end{equation}
where $\xi\in \Gamma(V^{\perp})$. Differentiating shows
\begin{equation}
\left\{\begin {array}{ll}\hat{Y}_{u}=b\hat{Y}+ \rho_{1}
(Y_u+bY)+\theta_{1}(Y_v+aY) +\eta_1+(\langle
\xi,\eta_1\rangle)Y,\\
 \hat{Y}_{v}=a\hat{Y}+ \theta_{2}
(Y_u+bY)+\rho_{2}(Y_v+aY) +\eta_2+(\langle \xi,\eta_2\rangle)Y.
\end {array}\right.\end{equation}
Here
\begin{equation}\label{eq-def}
\left\{\begin {array}{ll} \rho_{1}=2a_{u}-2\langle
\kappa_1,\kappa_2\rangle +\frac{1}{2}\langle
\xi,\xi\rangle,~~\rho_{2}=2b_{v}-2\langle \kappa_1,\kappa_2\rangle
+\frac{1}{2}\langle
\xi,\xi\rangle;\\
\theta_{1}=2b_{u}-2 b^{2} -s_1-2\langle
\xi,\kappa_1\rangle,~~\theta_{2}=2a_{v}-2a^{2}-s_2-2\langle
\xi,\kappa_2\rangle;\\
\eta_1=D_u\xi-b\xi+2D_v\kappa_1+2a\kappa_1,~~\eta_2=D_v\xi-a\xi+2D_u\kappa_2+2b\kappa_2.
\end {array}\right.
\end{equation}
By \eqref{eq-blas}, there must be $\eta_1=\eta_2=0$ and
$\rho_{1}=\rho_{2}=0$ or $\eta_1=\eta_2=0$ and
$\theta_{1}=\theta_{2}=0$.

From $\eta_1=\eta_2=0$, we obtain
$$D_v\kappa_1=-\frac{1}{2}D_u\xi+\frac{b}{2}\xi-a\kappa_1,D_u\kappa_2=-\frac{1}{2}D_v\xi+\frac{a}{2}\xi-b\kappa_2.$$
So
\begin{align*}
D_vD_v\kappa_1 +&\frac{s_2}{2}\kappa_1~=~D_v(-\frac{1}{2}D_u\xi+\frac{b}{2}\xi-a\kappa_1)+\frac{s_2}{2}\kappa_1\\
=&~-(\frac{\theta_2}{2}+2\langle \xi,\kappa_2\rangle)\kappa_1
+(\frac{b_v}{2}-\frac{ab}{2})\xi+\frac{a}{2}D_u\xi+\frac{b}{2}D_v\xi-\frac{1}{2}D_vD_u\xi.
\end{align*}
And
\begin{align*}
D_uD_u\kappa_2 +&\frac{s_1}{2}\kappa_2~=~D_u(-\frac{1}{2}D_v\xi+\frac{a}{2}\xi-b\kappa_2)+\frac{s_1}{2}\kappa_2\\
=&~-(\frac{\theta_1}{2}+2\langle
\xi,\kappa_1\rangle)\kappa_2 +(\frac{a_u}{2}-\frac{ab}{2}\xi+\frac{b}{2}D_v\xi+\frac{a}{2}D_u\xi-\frac{1}{2}D_uD_v\xi.\\
\end{align*}Plus the conformal Codazzi equation \eqref{eq-C2} and conformal Ricci
equation \eqref{eq-R2}, we get
\begin{equation}\label{eq-isoth-t}
\frac{a_u}{2}\xi-\frac{\theta_1}{2}\kappa_2=\frac{b_v}{2}\xi-\frac{\theta_2}{2}\kappa_1.
\end{equation}
This equation works when concerning the isothermic case.

Besides this, by the conformal Gauss equation \eqref{eq-G2}, we see
that
\begin{align*}
\theta_{1v}~=&~2b_{uv}-4 bb_{v} -s_{1v}-2\langle
\xi,\kappa_1\rangle_v\\~
~=&~\rho_{2u}-2b\rho_{2}+\langle-4D_u\kappa_2-4b\kappa_2-2D_v\xi
,\kappa_1\rangle-\langle
\xi,2D_v\kappa_1-b\xi+D_u\xi\rangle\\
~=&~\rho_{2u}-2b\rho_{2}.
\end{align*}
i.e.\begin{equation}\label{eq-theta-rho-1} \theta_{1v}=
\rho_{2u}-2b\rho_{2}.
\end{equation}
Similarly we obtain
\begin{equation}\label{eq-theta-rho-2}
\theta_{2u}= \rho_{1v}-2a\rho_{1}.
\end{equation}

Now let us prove Theorem A in the following three cases.\\

{\bf 1.\ \ The $S-Willmore$ case: $\theta_{1}=\theta_{2}=0,\
\xi=0$}\\

Since $\xi=0$ and $\eta_1=\eta_2=0$, \eqref{eq-def} reduces to
\begin{equation}
2D_v\kappa_1+a\kappa_1=2D_u\kappa_2+b\kappa_2=0.
\end{equation}
Together with
\begin{equation}\left\{\begin {array}{ll}
D_vD_v\kappa_1+\frac{s_2}{2}\kappa_1=\theta_{1}\kappa_1=0,\\
D_uD_u\kappa_2+\frac{s_1}{2}\kappa_2=\theta_{2}\kappa_2=0,\end
{array}\right.
\end{equation}
we see that $Y$ is a timelke S-Willmore surface. To verify
$\hat{Y}$, direct calculation shows that
$$\hat{Y}_{uv}=(\cdots)Y\mod \{\hat{Y},\hat{Y}_u,\hat{Y}_v\},\
 \hat{\kappa}_1=\rho_1\kappa_1,\ \hat{\kappa}_2=\rho_2\kappa_2.$$
 So $\hat{y}$ is S-Willmore by Theorem 2.2 since $\hat{Y}$ shares the same asymptotic coordinate and the same conformal Gauss map with $Y$.
So $y$ and $\hat{y}$
are a pair of dual S-Willmore surfaces.\\

{\bf 2.\ \ The isothermic case: $\rho_1=\rho_{2}=0$}\\

From the definition of $\rho_1$ and $\rho_{2}$, we see that
$a_u=b_v$. Substituting into \eqref{eq-isoth-t}
obtains\begin{equation} \theta_1\kappa_2=\theta_2\kappa_1.
\end{equation}
By use of \eqref{eq-theta-rho-1} and \eqref{eq-theta-rho-2}, we
have$$\theta_{1v}=\theta_{2u}=0.$$So$$\theta_{1}=\theta_{1}(u),\
\theta_{2}=\theta_{2}(v).$$ By choosing new asymptotic coordinate
$(\tilde{u},\tilde{v})$ we can derive
\begin{equation}
\tilde{\kappa}_1=\theta_2\kappa_1=\theta_1\kappa_2=\pm\tilde{\kappa}_2,
\end{equation}
where $\pm$ corresponds to $(\pm)-$isothemic surface. Notice that we must choose the $(\tilde{u},\tilde{v})$ such that
$$\left|\frac{\partial(\tilde{u},\tilde{v})}{\partial(u,v)}\right|>0$$
to ensure that $\langle Y_{\tilde{u}},Y_{\tilde{v}}\rangle>0$.

To show that $\hat{Y}$ is also $(\pm)-$isothemic surface as $y$, we
can suppose that $\kappa_1=\pm\kappa_2$. So
$\theta_{1}=\pm\theta_{2}$ and $\theta_{1v}=\theta_{2u}=0$ show that
$\theta_{1}=\pm\theta_{2}=\theta=const$. Then
$$\hat{Y}_u=b\hat{Y}+\theta(Y_v+aY)\Rightarrow Y_v=-aY+\frac{1}{\theta}(\hat{Y}_u-b\hat{Y}),$$
$$\hat{Y}_v=a\hat{Y}\pm\theta(Y_u+bY)\Rightarrow Y_u=-bY\pm\frac{1}{\theta}(\hat{Y}_v-a\hat{Y}).$$
So $\hat{Y}$ also satisfies the conditions of case 2, which means $\hat{Y}$ is also $(\pm)-$isothemic as $Y$.
In fact, $\hat{y}$ is the Darboux transform of $\theta-$parameter of $y$ and vice versa. \\

{\bf 3.\ \ The trivial case: $\theta_{1}=\theta_{2}=0, \
\xi\neq0$}\\

In this case, from \eqref{eq-theta-rho-1} and
\eqref{eq-theta-rho-2}, we can see that
$$\rho_{1v}-2a\rho_{1}=\rho_{2u}-2b\rho_{2}=0.$$
Together with \eqref{eq-isoth-t}, we see that $a_u=b_v$. So
$\rho_1=\rho_{2}=\rho\neq0$. Consider the vector
$\frac{1}{\rho}\hat{Y}-Y$, we have
\[
(\frac{1}{\rho}\hat{Y}-Y)_u~=~-b(\frac{1}{\rho}\hat{Y}-Y),
\
(\frac{1}{\rho}\hat{Y}-Y)_v~=~-a(\frac{1}{\rho}\hat{Y}-Y).
\]
This means that $\frac{1}{\rho}\hat{Y}-Y$ is a fixed direction,
showing that this is the trivial case.

\section{Proof of Theorem B}

Let $y:M\rightarrow Q^3_1$ be a timelike $(+)-$isothermic Willmore
surface with the adapted asymptotic coordinate $(u,v)$ and canonical
lift $Y$. Then $$\kappa_1=\kappa_2,
D_vD_v\kappa_1+\frac{s_2}{2}\kappa_1=D_uD_u\kappa_2+\frac{s_1}{2}\kappa_2=0.$$

 Assume that $$\kappa_1=\kappa_2=kE,$$
where $E$ is a unit section of the conformal normal bundle. If $k=0$
in a neighborhood, $y$ is contained in some $S^2_1$ and then minimal
in some $S^3_1$.

So we can suppose $k\neq0$ in a open subset $U\subset M$. Set
\begin{equation}\hat{Y}=N+2aY_u+2bY_v+(2ab)Y,
\end{equation}
with $$a=-\frac{\kappa_v}{\kappa},\ b=-\frac{\kappa_u}{\kappa}.$$
From the calculation in Section 3, we can verify that $\hat{Y}$ is
just the dual Willmore surface of $Y$ and
$$\hat{Y}_u=a\hat{Y}+\rho(Y_u+bY),\hat{Y}_v=b\hat{Y}+\rho(Y_v+aY),$$
where $$\rho=a_u-2k^2=b_v-2k^2,\ \rho_u=2b\rho, \rho_v=2a\rho,$$ by
use of the Willmore equations as above in Section 3.

Consider the vector field $$Y_0=\hat{Y}-\rho Y. $$ Differentiating
it leads to
$$Y_{0u}=\hat{Y}_u-\rho_uY-\rho Y_u=b(\hat{Y}-\rho Y)=bY_0,\ Y_{0v}=\hat{Y}_v-\rho_vY-\rho Y_v=aY_0.$$
This means that $Y_0$ is a point in $Q^3_1$.

(i) If $\langle Y_0,Y_0\rangle=0$, $\rho\equiv0$. So $[\hat{Y}]$
reduces to a point. By some conformal transform, we can set
$\hat{Y}=f_1(1,0,0,0,1)$ with some function $f_1$ and
$Y=e^{-\omega}(\frac{-1+\langle x,x\rangle}{2},x,\frac{1+\langle
x,x\rangle}{2})$ for some timelike surface $x:U\rightarrow R^3_1$
with $\langle x_u,x_v\rangle=\frac{1}{2}e^{2\omega}$.

The structure equations of $x$ is:
\begin{equation*}
\left\{ \begin{aligned}
         x_{uu} &= 2\omega_ux_u+\Omega n,  x_{vv} = 2\omega_vx_v+\Omega n,
                 x_{uv} =\frac{1}{2}e^{2\omega}Hn,\\
                 n_{u}&=-Hx_u-2\Omega e^{-2\omega}x_v,  n_{v}=-2\Omega
                 e^{-2\omega}x_u-Hx_v.
                          \end{aligned} \right.
                          \end{equation*}
So
$$k=e^{-\omega}\Omega, a=\omega_v-\frac{\Omega_v}{\Omega},b=\omega_u-\frac{\Omega_u}{\Omega},$$
$$N=e^\omega(1+H\langle x,n\rangle,Hn,1+H\langle x,n\rangle)-2\omega_vY_u-2\omega_uY_v+2\omega_u\omega_vY,$$
$$\hat{Y}=e^\omega(1+H\langle x,n\rangle,Hn,1+H\langle x,n\rangle)+(\cdots)Y_u+(\cdots)Y_v+(\cdots)Y.$$
Since $\hat{Y}=f_1(1,0,0,0,1)$, the coefficient of $n$ must be zero,
i.e. $H=0$, which means that $x$ is a timelike minimal surface in
$R^3_1$.

(ii) If $\langle Y_0,Y_0\rangle<0$, by some conformal transform, we
can set $Y_0=f_2(0,0,0,0,1)$ with some function $f_2$ and
$Y=e^{-\omega}(x,1)$ for some timelike surface $x:U\rightarrow
S^3_1$ with $\langle x_u,x_v\rangle=\frac{1}{2}e^{2\omega}$. Similar
to case (i), it is direct to show that $x$ is just a minimal surface
in $S^3_1$.

(iii) If $\langle Y_0,Y_0\rangle>0$, similar treatments as above
show that $Y_0=f_3(1,0,0,0,0)$ with some function $f_3$ and
$Y=e^{-\omega}(1,x)$ for some timelike minimal surface
$x:U\rightarrow H^3_1$ with $\langle x_u,x_v\rangle=\frac{1}{2}e^{2\omega}$.\\\\\\ {\bf Acknowledgement}\\\\
The author would like to thank Professor Changping Wang for his
direction and help. It is also a pleasure to thank Dr. Xiang Ma for
his suggestions and discussions.

\def\refname{Reference}

\end{document}